\documentclass[11pt]{amsart}
\usepackage[T1]{fontenc}
\usepackage[utf8]{inputenc}
\usepackage{lmodern}
\usepackage{amsmath,amssymb,mathtools}
\usepackage[letterpaper,textwidth=360pt,textheight=584pt,centering]{geometry}
\usepackage{microtype}

\usepackage{hyperref}
\hypersetup{hypertex=true,
            colorlinks=true,
            linkcolor=blue,
            anchorcolor=blue,
            citecolor=blue }

\numberwithin{equation}{section}
\newtheorem{theorem}{Theorem}[section]
\newtheorem{lemma}[theorem]{Lemma}
\newtheorem{proposition}[theorem]{Proposition}
\theoremstyle{remark}
\newtheorem{remark}[theorem]{Remark}
\newcommand{\R}{\mathbb R}
\newcommand{\EL}{\mathcal E_{L,p}}
\newcommand{\Sp}{\mathsf S_p}
\newcommand{\dd}{\,\mathrm d}
\newcommand{\norm}[1]{\left\lVert #1\right\rVert}
\title[Strict logarithmic rearrangement inequalities]{Strict rearrangement inequalities for logarithmic $p$-Laplacians}

\author{Hichem Hajaiej}
\address{Hichem Hajaiej: Department of Mathematics, California State University, Los Angeles, USA}
\email{hichem.hajaiej@gmail.com}

\author{Qinyang Li}
\address{Qinyang Li: School of Mathematical Sciences, Xiamen University, Xiamen 361005, P. R. China}
\email{19020241153813@stu.xmu.edu.cn}

\author{Weicheng Zhan}
\address{Weicheng Zhan: School of Mathematical Sciences, Xiamen University, Xiamen, Fujian 361005, P. R. China}
\email{zhanweicheng@amss.ac.cn}

\begin{document}
\begin{abstract}
We prove a P\'olya--Szeg\H{o} inequality for the logarithmic $p$-Laplacian, $1<p<\infty$, for real-valued functions in its natural form space. The inequality holds under a logarithmic integrability condition, without any restriction on the support. This condition is also necessary for the rearranged function to have finite logarithmic energy. For a nonzero function, equality holds if and only if the function has one sign and its modulus is radially nonincreasing about some point. As a consequence, balls uniquely minimize the first Dirichlet eigenvalue among bounded open sets of prescribed measure, up to translations and null sets. These results answer the equality questions raised by Chen and Weth~\cite{CW} for $p=2$ and by Dyda, Jarohs and Sk~\cite{DJS} for general $p$. The proof combines strict rearrangement for supermodular interactions with monotonicity of regularized interaction deficits, which allows us to retain the equality information in the singular limit.


\bigskip

\end{abstract}
\maketitle

\section{Introduction}\label{s1}
\subsection{Background and main ideas}
The logarithmic Laplacian $L_\Delta$ arises by differentiating the fractional Laplacian with respect to its order at zero. For $u\in C_c^\infty(\R^N)$, we have
$$
(-\Delta)^s u=u+sL_\Delta u+o(s)\qquad\text{as }s\to0^+,
$$
and the Fourier symbol of $L_\Delta$ is $2\log|\xi|$. Chen and Weth~\cite{CW} developed the variational and spectral theory of the corresponding Dirichlet problem on bounded domains. This theory was extended to the nonlinear setting by Dyda, Jarohs and Sk~\cite{DJS}, who introduced the logarithmic $p$-Laplacian for $1<p<\infty$.

Several aspects of the nonlinear theory have since been studied. Arora, Giacomoni, Hajaiej and Vaishnavi~\cite{AGHV} established the $p$-logarithmic Sobolev inequality and sharp embeddings into Orlicz-type spaces, and studied Dirichlet problems with critical growth. Arora, Hajaiej and Perera~\cite{AHP} constructed an unbounded sequence of minimax eigenvalues and obtained existence results for superlinear subcritical Dirichlet problems. For positive fractional orders, Bahrouni, Gouasmia, Hajaiej and Ouannasser~\cite{BGHO} studied the fractional logarithmic $p$-Laplacian, obtained by differentiating the fractional $p$-Laplacian with respect to its order at $s>0$.

Symmetrization for the logarithmic energy requires some care. The energy consists of a positive short-range term, a renormalized long-range interaction, and a multiple of the $L^p$-mass. Its kernel $|x-y|^{-N}$ is singular at the diagonal and is not integrable at infinity. Moreover, for real-valued functions the numerator of the long-range interaction may have either sign. The cancellation in the energy must therefore be preserved, and rearrangement cannot be applied separately to terms that are only formally defined.

Chen and Weth~\cite[Corollary 3.6]{CW} proved that balls minimize the first Dirichlet eigenvalue among bounded Lipschitz domains of prescribed measure. Their argument uses a first-order expansion of fractional Dirichlet eigenvalues. They left open the equality case and the question of a direct symmetrization proof; see~\cite[Remark 3.7]{CW}. By a similar asymptotic argument, Dyda, Jarohs and Sk obtained the Faber--Krahn inequality for the logarithmic $p$-Laplacian and asked whether equality characterizes balls~\cite[Remark 1.5]{DJS}. Although fractional eigenvalue asymptotics imply the spectral comparison, the equality information need not pass to the limit. Indeed, a fractional spectral deficit that is positive for every $s>0$ may still be $o(s)$ as $s\to0^+$. Thus the equality results for fractional P\'olya--Szeg\H{o} inequalities, such as~\cite[Theorem A.1]{FS}, do not by themselves resolve these questions.

In this paper, we establish a direct P\'olya--Szeg\H{o} inequality for the logarithmic $p$-energy and determine all equality cases for real-valued functions. A nonzero function attains equality precisely when it has one sign and its modulus is a translate of its Schwarz rearrangement. The result holds on $\R^N$ without any boundedness or finite-measure assumption on the support. Applying it to a positive first eigenfunction, we prove that balls are the unique minimizers of the first Dirichlet eigenvalue among bounded open sets of fixed measure, up to translations and null sets. In particular, no boundary regularity is needed for this conclusion.

An essential point is to identify the integrability condition under which the rearranged energy is finite. We assume that
$$
\Sp(u):=\int_{\R^N}|u|^p\log^-|u|\dd x<\infty,\qquad \log^-t:=\max\{0,-\log t\}.
$$
This condition controls small values of the function on sets that may be unbounded. It is invariant under rearrangement and is weaker than the additional assumption $u\in L^q(\R^N)$ for some $0<q<p$. We prove that, for $u\in X^p(\R^N)$, it is necessary and sufficient for $u^\#=(|u|)^*$ to have finite logarithmic energy. More precisely, it is equivalent both to a logarithmic spatial moment of $u^\#$ and to finiteness of its long-range interaction. If the condition fails, the short-range term of $u^\#$ remains finite, whereas its long-range term equals $-\infty$.

To prove the inequality and its equality statement, we write the interaction in terms of the strictly supermodular function
$$
H_p(a,b)=a^p+b^p-|a-b|^p,\qquad a,b\ge0,
$$
and regularize the kernel by $K_\varepsilon(z)=(|z|^2+\varepsilon^2)^{-N/2}$. We first estimate these interactions in terms of the distribution of the function and $\Sp(u)$. The logarithmic energy deficit can then be expressed as the limit of the regularized rearrangement deficits. The key observation is that these deficits increase as $\varepsilon\downarrow0$, because the difference of two regularized kernels is positive and strictly radially decreasing. Consequently, equality in the limiting comparison forces equality at every fixed regularization scale. By strict rearrangement at one such scale, the function is a translate of its Schwarz rearrangement. A separate comparison between $u$ and $|u|$ determines the sign of a real-valued equality case.

The use of monotonicity in this argument is motivated by the work of De Luca, Novaga and Ponsiglione~\cite{DNP}. For characteristic functions in the linear case, they proved that balls uniquely minimize the $0$-fractional perimeter at fixed volume; see~\cite[Theorem 7.2]{DNP}. Their proof also uses monotonicity of regularized interaction differences. Here we combine this property with the strict rearrangement theorem of Burchard and Hajaiej~\cite{BH} for supermodular integrands. This allows us to treat general real-valued functions, all $1<p<\infty$, and the precise logarithmic finiteness condition. For the classical equality theory of the Riesz rearrangement inequality, we refer to Burchard~\cite{B}.

\subsection{Main results}
Let $N\ge1$ and $1<p<\infty$. Following~\cite{DJS}, we set
\begin{equation}\label{1-1}
C_{N,p}:=\frac{p\Gamma(N/2)}{2\pi^{N/2}}=\frac p{|S^{N-1}|}, \qquad \rho_{N,p}:=2\log2-\gamma+\frac p2\psi\!\left(\frac N2\right),
\end{equation}
where $\gamma$ is Euler's constant and $\psi=\Gamma'/\Gamma$. We write
\begin{equation}\label{1-2}
[u]_{X^p}^p:=\iint_{|x-y|<1}\frac{|u(x)-u(y)|^p}{|x-y|^N}\dd x\dd y,
\end{equation}
and define $X^p(\R^N):=\{u\in L^p(\R^N):[u]_{X^p}<\infty\}$. For a nonnegative function $u$, let $u^*$ denote its Schwarz rearrangement. For real-valued $u$, set $u^\#:=(|u|)^*$.

We also set
\begin{equation}\label{1-3}
\Sp(u):=\int_{\R^N}|u(x)|^p\log^-|u(x)|\dd x,
\end{equation}
with the convention $0^p\log^-0=0$. This quantity is unchanged by rearrangement. Its finiteness is also preserved under multiplication by a nonzero constant, since
$$
|\log^-(ct)-\log^-t|\le|\log c|\qquad(c,t>0).
$$
For $u\in X^p(\R^N)$ satisfying $\Sp(u)<\infty$, the logarithmic $p$-energy is
\begin{equation}\label{1-4}
\begin{aligned}
\EL(u):={}&\frac{C_{N,p}}2\iint_{|x-y|<1}\frac{|u(x)-u(y)|^p}{|x-y|^N}\dd x\dd y\\
&+\frac{C_{N,p}}2\iint_{|x-y|\ge1}\frac{|u(x)-u(y)|^p-|u(x)|^p-|u(y)|^p}{|x-y|^N}\dd x\dd y\\
&+\rho_{N,p}\int_{\R^N}|u|^p\dd x.
\end{aligned}
\end{equation}
The absolute convergence of the long-range integral follows from Lemma~\ref{lem2-2} and~\eqref{2-3} below. With this definition, Schwarz rearrangement decreases the energy, and equality determines both the sign and the radial profile of the function.

\begin{theorem}\label{thm1-1}
Let $N\ge1$ and $1<p<\infty$. Suppose that $u\in X^p(\R^N)$ is real-valued and satisfies $\Sp(u)<\infty$. Then $u^\#\in X^p(\R^N)$, $\Sp(u^\#)=\Sp(u)$, and
\begin{equation}\label{1-5}
\EL(u^\#)\le\EL(u).
\end{equation}
If $u\not\equiv0$, equality holds in~\eqref{1-5} if and only if there exist $x_0\in\R^N$ and $\sigma\in\{-1,1\}$ such that
$$
u=\sigma u^\#(\cdot-x_0)\qquad\text{a.e. in }\R^N.
$$
\end{theorem}

The logarithmic condition is satisfied whenever $|\{u\ne0\}|<\infty$, even if this set is unbounded. It is also satisfied if $u\in L^q(\R^N)$ for some $0<q<p$. Indeed, for $0<t<1$,
\begin{equation}\label{1-6}
t^p\log\frac1t\le\frac1{pe},\qquad t^p\log\frac1t\le\frac{t^q}{e(p-q)}.
\end{equation}
Hence Theorem~\ref{thm1-1} applies to all compactly supported functions in $X^p(\R^N)$ and to all functions in $X^p(\R^N)\cap L^1(\R^N)$. Remark~\ref{rem2-4} contains admissible functions that satisfy none of these stronger assumptions.

The logarithmic condition is also necessary for the rearranged energy to be finite. More precisely, for every real-valued $u\in X^p(\R^N)$, one has $u^\#\in X^p(\R^N)$ and
\begin{equation}\label{1-7}
\EL(u^\#)\in\R\quad\Longleftrightarrow\quad\Sp(u)<\infty.
\end{equation}
If $\Sp(u)=\infty$, the short-range term of $u^\#$ remains finite, whereas its long-range term equals $-\infty$. These assertions follow from Proposition~\ref{prop2-3} and Lemma~\ref{lem3-2}.

We next consider the first Dirichlet eigenvalue. For a bounded nonempty open set $\Omega\subset\R^N$, let
\begin{equation}\label{1-8}
X_0^p(\Omega):=\{u\in X^p(\R^N):u=0\text{ a.e. in }\R^N\setminus\Omega\}.
\end{equation}
The first Dirichlet eigenvalue is defined by
\begin{equation}\label{1-9}
\lambda_{1,L,p}(\Omega):=\inf_{u\in X_0^p(\Omega)\setminus\{0\}} \frac{\EL(u)}{\norm{u}_{L^p(\R^N)}^p}.
\end{equation}
By~\eqref{1-6}, every function in $X_0^p(\Omega)$ satisfies the logarithmic integrability condition. Theorem~\ref{thm1-1} therefore yields the following Faber--Krahn inequality and equality classification.

\begin{theorem}\label{thm1-2}
Let $N\ge1$, $1<p<\infty$, and let $\Omega\subset\R^N$ be a bounded nonempty open set. If $B\subset\R^N$ is a ball with $|B|=|\Omega|$, then
\begin{equation}\label{1-10}
\lambda_{1,L,p}(B)\le\lambda_{1,L,p}(\Omega).
\end{equation}
Equality holds if and only if $|\Omega\mathbin\triangle(B+x_0)|=0$ for some $x_0\in\R^N$. In particular, if $\Omega$ has Lipschitz boundary, equality holds if and only if $\Omega$ is a ball.
\end{theorem}

For general open sets, equality of the domains must be understood up to null sets: both the space in~\eqref{1-8} and the quotient in~\eqref{1-9} are unchanged by null modifications of the domain. The Lipschitz assumption in the last assertion allows us to recover equality of the domains themselves.

The equality case in Theorem~\ref{thm1-1} rests on a comparison of regularized interaction deficits. For $u\ge0$, define
$$
J_\varepsilon(u):=\iint_{\R^N\times\R^N} K_\varepsilon(x-y)H_p(u(x),u(y))\dd x\dd y.
$$
Under the hypotheses of Theorem~\ref{thm1-1}, these interactions are finite for both $u$ and $u^*$. We prove that
\begin{equation}\label{1-11}
\EL(u)-\EL(u^*)=\frac{C_{N,p}}2\lim_{\varepsilon\downarrow0} \bigl(J_\varepsilon(u^*)-J_\varepsilon(u)\bigr).
\end{equation}
The deficits on the right-hand side are nonnegative and increase as $\varepsilon\downarrow0$. Hence the energy deficit bounds the Riesz deficit at each fixed scale from below. When the energy deficit vanishes, we may apply the strict rearrangement theorem at any one of these scales.

The paper is organized as follows. Section~\ref{s2} establishes the long-range estimates and characterizes finiteness of the long-range term after rearrangement. Section~\ref{s3} proves the regularized rearrangement inequality and preservation of $X^p(\R^N)$. Sections~\ref{s4} and~\ref{s5} prove Theorems~\ref{thm1-1} and~\ref{thm1-2}, respectively. Section~\ref{s6} discusses sharpness of the logarithmic condition, the equality case at $p=1$, and the geometric and scaling conclusions of Theorem~\ref{thm1-2}. Appendix~\ref{appA} presents an alternative proof of Theorem~\ref{thm1-1} by polarization.

\section{Long-range integrability}\label{s2}
We first establish the long-range estimates needed to define the energy in~\eqref{1-4}. Since $|z|^{-N}$ is not integrable at infinity, the logarithmic condition is used to control the contribution of small values of the function. Throughout this section, $C$ denotes a positive constant that may vary from line to line; its dependence is specified in each estimate.

For $a,b\ge0$, let
\begin{equation}\label{2-1}
H_p(a,b):=a^p+b^p-|a-b|^p,\qquad G(z):=\min\{1,|z|^{-N}\},
\end{equation}
where $G(0)=1$. Then $0\le H_p(a,b)\le a^p+b^p$ and $H_p(a,0)=0$. For $0\le b\le a$, we also have
\begin{equation}\label{2-2}
a^{p-1}b\le H_p(a,b)\le(p+1)a^{p-1}b.
\end{equation}
The lower bound follows from $(1-t)^p\le1-t$ for $0\le t\le1$. The upper bound follows from the mean value theorem and $b^p\le a^{p-1}b$.

To treat real-valued functions, we shall also use the following consequence of~\eqref{2-2}: for all $a,b\in\R$,
\begin{equation}\label{2-3}
\bigl||a-b|^p-|a|^p-|b|^p\bigr|\le C_pH_p(|a|,|b|).
\end{equation}
If $ab\ge0$, the left-hand side equals $H_p(|a|,|b|)$. If $ab<0$, set $A=\max\{|a|,|b|\}$ and $B=\min\{|a|,|b|\}$. By the mean value theorem,
$$
0\le(A+B)^p-A^p-B^p\le p2^{p-1}A^{p-1}B\le p2^{p-1}H_p(A,B).
$$
Thus~\eqref{2-3} holds with $C_p=\max\{1,p2^{p-1}\}$.

We begin with an estimate on sets of finite measure. Its logarithmic dependence on the measure will be used in the superlevel-set decomposition below.
\begin{lemma}\label{lem2-1}
There exists $C=C(N)>0$ such that, for every measurable set $E\subset\R^N$ of finite measure,
\begin{equation}\label{2-4}
\sup_{x\in\R^N}\int_EG(x-y)\dd y\le C\log(e+|E|).
\end{equation}
\end{lemma}
\begin{proof}
Set $R=\max\{1,|E|^{1/N}\}$. By the definition of $G$,
$$
\int_EG(x-y)\dd y\le\int_{B_R}G(z)\dd z+R^{-N}|E| \le|B_1|+|S^{N-1}|\log R+1.
$$
The right-hand side is bounded by $C\log(e+|E|)$, independently of $x$. This proves~\eqref{2-4}.
\end{proof}

We next estimate the interaction against $G$ in terms of $\Sp(v)$. Since the bound depends only on the distribution of $v$, it remains valid after rearrangement.
\begin{lemma}\label{lem2-2}
Let $v\in L^p(\R^N)$ be nonnegative and satisfy $\Sp(v)<\infty$. Then
\begin{equation}\label{2-5}
J_G(v):=\iint_{\R^N\times\R^N}G(x-y)H_p(v(x),v(y))\dd x\dd y<\infty.
\end{equation}
Moreover, if $\norm{v}_{L^p(\R^N)}=1$, then
\begin{equation}\label{2-6}
J_G(v)\le C\bigl(1+\Sp(v)\bigr),
\end{equation}
where $C=C(N,p)>0$.
\end{lemma}
\begin{proof}
The assertion is immediate if $v=0$. By the homogeneity of $H_p$ and the estimate following~\eqref{1-3}, we may assume that $\norm{v}_{L^p(\R^N)}=1$. Let $\mu(t):=|\{v>t\}|$. By Chebyshev's inequality, we have $\mu(t)\le t^{-p}$ for $t>0$.

For $a>0$ and $x\in\R^N$, Tonelli's theorem and Lemma~\ref{lem2-1} yield
\begin{equation}\label{2-7}
\begin{aligned}
\int_{\{0<v\le a\}}G(x-y)v(y)\dd y &\le\int_0^a\int_{\{v>t\}}G(x-y)\dd y\dd t\\
&\le C\int_0^a\log(e+\mu(t))\dd t\\
&\le C\int_0^a(1+\log^-t)\dd t\\
&\le Ca(1+\log^-a).
\end{aligned}
\end{equation}
The last inequality follows by direct integration, treating $(0,1)$ and $(1,a)$ separately when $a>1$.

By symmetry,~\eqref{2-2}, and~\eqref{2-7} with $a=v(x)$, we have
\begin{align*}
J_G(v)&\le2(p+1)\int_{\R^N}v(x)^{p-1} \int_{\{0<v\le v(x)\}}G(x-y)v(y)\dd y\dd x\\
&\le C\int_{\R^N}v(x)^p(1+\log^-v(x))\dd x\\
&=C\bigl(1+\Sp(v)\bigr).
\end{align*}
This proves~\eqref{2-6} and completes the proof.
\end{proof}

By Lemma~\ref{lem2-2} and~\eqref{2-3}, the long-range integral in~\eqref{1-4} is absolutely convergent for $u\in X^p(\R^N)$ with $\Sp(u)<\infty$. We shall also use a translation formula for the energy. For $w\in L^p(\R^N)$, set
\begin{equation}\label{2-8}
M_w:=\norm{w}_{L^p(\R^N)}^p,\qquad \tau_w(z):=\int_{\R^N}|w(x+z)-w(x)|^p\dd x.
\end{equation}
Then $0\le\tau_w(z)\le2^pM_w$. For nonnegative $w$,
\begin{equation}\label{2-9}
I_w(z):=\int_{\R^N}H_p(w(x),w(x+z))\dd x =2M_w-\tau_w(z)\ge0.
\end{equation}
Let $w\in X^p(\R^N)$ be real-valued and satisfy $\Sp(w)<\infty$. By the preceding integrability estimates, we may make the change of variables $y=x+z$ in~\eqref{1-4} and obtain
\begin{equation}\label{2-10}
\begin{aligned}
\EL(w)={}&\frac{C_{N,p}}2\int_{|z|<1}\frac{\tau_w(z)}{|z|^N}\dd z\\
&+\frac{C_{N,p}}2\int_{|z|\ge1}\frac{\tau_w(z)-2M_w}{|z|^N}\dd z +\rho_{N,p}M_w.
\end{aligned}
\end{equation}
In particular, the energy is invariant under translations and multiplication by $-1$.

For rearranged functions, the logarithmic condition is also necessary. The following result requires only $L^p$-integrability.
\begin{proposition}\label{prop2-3}
Let $u\in L^p(\R^N)$ be real-valued, and put $v=u^\#$. Define
\begin{equation}\label{2-11}
T_p(v):=\iint_{|x-y|\ge1}\frac{H_p(v(x),v(y))}{|x-y|^N}\dd x\dd y\in[0,\infty].
\end{equation}
The following three conditions are equivalent:
\begin{equation}\label{2-12}
\Sp(u)<\infty,\qquad \int_{\R^N}v(x)^p\log(e+|x|)\dd x<\infty,\qquad T_p(v)<\infty.
\end{equation}
\end{proposition}
\begin{proof}
By equimeasurability, $\Sp(v)=\Sp(u)$. The bound $T_p(v)\le J_G(v)$ and Lemma~\ref{lem2-2} show that the first condition implies the third.

Suppose that the third condition holds. Choose a radial nonincreasing representative of $v$. For $x\in\R^N$ with $r=|x|\ge4$, let
$$
A_x:=\left\{z\in\R^N:1\le|z|\le\frac r2,\quad x\cdot z\le-\frac12r|z|\right\}.
$$
For $z\in A_x$, one has $|x+z|\le|x|$, hence $v(x+z)\ge v(x)$. It follows that
$$
H_p(v(x),v(x+z))\ge v(x)^p.
$$
The set of directions defining $A_x$ has positive measure depending only on $N$; when $N=1$, it consists of the direction towards the origin. Consequently,
$$
\int_{A_x}\frac{\dd z}{|z|^N}=c_N\log(r/2)\ge\frac{c_N}2\log r.
$$
By Tonelli's theorem, we have
\begin{equation}\label{2-13}
T_p(v)\ge c\int_{|x|\ge4}v(x)^p\log|x|\dd x,
\end{equation}
where $c=c(N)>0$. Since $v\in L^p(\R^N)$, the logarithmically weighted integral over $|x|<4$ is finite. Thus the second condition holds.

It remains to show that the second condition implies the first. For $t,r\ge0$,
\begin{equation}\label{2-14}
t^p\log^-t\le\frac{N+1}{p}t^p\log(e+r)+\frac1{pe}(e+r)^{-N-1}.
\end{equation}
For $t\ge1$, this is immediate. For $0<t<1$, set $a=t^p$ and $b=(e+r)^{-N-1}$. Since $s\log(1/s)\le1/e$ for $s>0$,
$$
a\log\frac1a=a\log\frac1b+b\frac ab\log\frac ba \le(N+1)a\log(e+r)+\frac be.
$$
After division by $p$, this is~\eqref{2-14}. The inequality also holds at $t=0$ by continuity. Applying it with $t=v(x)$ and $r=|x|$ and integrating over $\R^N$, we obtain $\Sp(v)<\infty$.
\end{proof}

The next example shows both that the logarithmic condition cannot be omitted and that it is weaker than additional $L^q$-integrability for $q<p$.
\begin{remark}\label{rem2-4}
For $\alpha>1$, set
\begin{equation}\label{2-15}
v_\alpha(x):=\frac1{(1+|x|^2)^{N/(2p)}[\log(e+|x|^2)]^{\alpha/p}}.
\end{equation}
Then $v_\alpha$ is positive, smooth and radially nonincreasing. As $r=|x|\to\infty$,
$$
v_\alpha(x)^p\asymp r^{-N}(\log r)^{-\alpha},\qquad |\nabla v_\alpha(x)|^p\le Cr^{-N-p}(\log r)^{-\alpha}.
$$
Thus $v_\alpha\in W^{1,p}(\R^N)\subset X^p(\R^N)$. The latter inclusion follows by integrating $\tau_w(z)\le|z|^p\norm{\nabla w}_{L^p(\R^N)}^p$ over $B_1$ against $|z|^{-N}$. Moreover, $v_\alpha\notin L^q(\R^N)$ for every $0<q<p$, whereas
$$
\Sp(v_\alpha)<\infty\quad\Longleftrightarrow\quad \int^\infty\frac{\dd r}{r(\log r)^{\alpha-1}}<\infty \quad\Longleftrightarrow\quad\alpha>2.
$$
Hence Theorem~\ref{thm1-1} applies when $\alpha>2$. If $1<\alpha\le2$, then $T_p(v_\alpha)=\infty$ by Proposition~\ref{prop2-3}, and the expression in~\eqref{1-4} equals $-\infty$.
\end{remark}

\section{A regularized rearrangement inequality}\label{s3}
We now apply the strict rearrangement theorem for supermodular integrands to $H_p$. The resulting comparison will also be used with truncated kernels to prove that Schwarz rearrangement preserves $X^p(\R^N)$.

The interaction $H_p$ in~\eqref{2-1} is strictly supermodular. More precisely, for $a,b\ge0$ and $h,k>0$,
\begin{equation}\label{3-1}
\begin{aligned}
&H_p(a+h,b+k)+H_p(a,b)-H_p(a+h,b)-H_p(a,b+k)\\
&\hspace{2em}=p(p-1)\int_0^h\int_0^k|a-b+s-t|^{p-2}\dd t\dd s>0.
\end{aligned}
\end{equation}
Indeed, the derivative of $t\mapsto|t|^p$ is locally absolutely continuous for every $p>1$. Two applications of the fundamental theorem of calculus imply~\eqref{3-1}. When $1<p<2$, the possible singularity at $a-b+s-t=0$ is locally integrable, so the argument remains valid.

For a bounded nonnegative nonincreasing function $\kappa:[0,\infty)\to[0,\infty)$ and a nonnegative $w\in L^p(\R^N)$, define
\begin{equation}\label{3-2}
J_\kappa(w):=\iint_{\R^N\times\R^N}\kappa(|x-y|)H_p(w(x),w(y))\dd x\dd y\in[0,\infty].
\end{equation}
The following proposition is the form of the strict rearrangement theorem of Burchard and Hajaiej~\cite[Theorem 2]{BH} needed here.
\begin{proposition}\label{prop3-1}
Let $w\in L^p(\R^N)$ be nonnegative, and let $\kappa:[0,\infty)\to[0,\infty)$ be bounded and nonincreasing. Then
\begin{equation}\label{3-3}
J_\kappa(w)\le J_\kappa(w^*).
\end{equation}
Suppose, in addition, that $w\not\equiv0$, that $\kappa$ is positive and strictly decreasing, and that $J_\kappa(w^*)<\infty$. Equality holds in~\eqref{3-3} if and only if $w=w^*(\cdot-x_0)$ almost everywhere for some $x_0\in\R^N$.
\end{proposition}
\begin{proof}
Since $w\in L^p(\R^N)$, all its positive superlevel sets have finite measure. By~\eqref{3-1}, $H_p$ is strictly supermodular; it is also continuous, nonnegative and vanishes when either variable is zero. We may therefore apply~\cite[Theorem 2]{BH} with $m=2$, $u_1=u_2=w$, $F=H_p$, and $K_{12}=\kappa$ to obtain~\eqref{3-3}.

Under the additional assumptions, both interactions are finite, and the nonzero $L^p$ function $w$ cannot be constant almost everywhere on $\R^N$. If equality holds in~\eqref{3-3}, the equality assertion of the same theorem implies that $w$ is a translate of $w^*$. The converse follows from translation invariance.
\end{proof}

Applying Proposition~\ref{prop3-1} to integrable truncated kernels, we obtain the next comparison of the short-range terms. This step requires no logarithmic integrability assumption.
\begin{lemma}\label{lem3-2}
Let $w\in X^p(\R^N)$ be real-valued. Then $w^\#\in X^p(\R^N)$ and
\begin{equation}\label{3-4}
[w^\#]_{X^p}^p\le[|w|]_{X^p}^p\le[w]_{X^p}^p.
\end{equation}
\end{lemma}
\begin{proof}
The second inequality follows from the reverse triangle inequality. We may therefore assume that $w\ge0$. For $0<\delta<1$, let
$$
k_\delta(z):=\mathbf1_{\{|z|<1\}}\min\{\delta^{-N},|z|^{-N}\}.
$$
The kernel $k_\delta$ is bounded, radially nonincreasing and belongs to $L^1(\R^N)$. Since $H_p(a,b)\le a^p+b^p$, both interactions against $k_\delta$ are bounded by $2M_w\norm{k_\delta}_1$. By Proposition~\ref{prop3-1} and~\eqref{2-9}, we have
$$
\int_{\R^N}k_\delta(z)(2M_w-\tau_w(z))\dd z \le\int_{\R^N}k_\delta(z)(2M_w-\tau_{w^*}(z))\dd z,
$$
where $M_{w^*}=M_w$. The common mass term is finite, so
$$
\int_{\R^N}k_\delta(z)\tau_{w^*}(z)\dd z \le\int_{\R^N}k_\delta(z)\tau_w(z)\dd z.
$$
Letting $\delta\downarrow0$ and applying the monotone convergence theorem proves~\eqref{3-4}.
\end{proof}

For $u\in X^p(\R^N)$, the short-range term of $u^\#$ is finite by Lemma~\ref{lem3-2}, whereas its long-range term is $-C_{N,p}T_p(u^\#)/2$. Proposition~\ref{prop2-3} therefore proves~\eqref{1-7}.

For $\varepsilon>0$, set
\begin{equation}\label{3-5}
\kappa_\varepsilon(r):=(r^2+\varepsilon^2)^{-N/2},\qquad K_\varepsilon(z):=\kappa_\varepsilon(|z|),\qquad J_\varepsilon:=J_{\kappa_\varepsilon}.
\end{equation}
By definition,
\begin{equation}\label{3-6}
K_\varepsilon(z)\le\max\{1,\varepsilon^{-N}\}G(z).
\end{equation}
Lemma~\ref{lem2-2} and equimeasurability imply that $J_\varepsilon(w)$ and $J_\varepsilon(w^*)$ are finite whenever $w\in L^p(\R^N)$ is nonnegative and $\Sp(w)<\infty$. By Tonelli's theorem and~\eqref{2-9}, we have
\begin{equation}\label{3-7}
J_\varepsilon(w)=\int_{\R^N}K_\varepsilon(z)I_w(z)\dd z =\int_{\R^N}K_\varepsilon(z)(2M_w-\tau_w(z))\dd z.
\end{equation}
The cancellation in $2M_w-\tau_w(z)$ is needed in this identity, since the two terms need not be integrable separately.

\section{Proof of Theorem~\ref{thm1-1}}\label{s4}
We first treat nonnegative functions. Using the estimates of Section~\ref{s2}, we express the energy difference as the limit of regularized interaction deficits. Their monotonicity will allow us to pass the equality information to the limit.

For nonnegative $u\in X^p(\R^N)$ with $\Sp(u)<\infty$, define
\begin{equation}\label{4-1}
D_\varepsilon(u):=J_\varepsilon(u^*)-J_\varepsilon(u),\qquad\varepsilon>0.
\end{equation}
By Lemma~\ref{lem2-2} and~\eqref{3-6}, both interactions are finite, and $D_\varepsilon(u)\ge0$ by Proposition~\ref{prop3-1}. The next lemma identifies the limit of these deficits and proves their monotonicity.
\begin{lemma}\label{lem4-1}
Let $u\in X^p(\R^N)$ be nonnegative and satisfy $\Sp(u)<\infty$. Then
\begin{equation}\label{4-2}
\begin{aligned}
\EL(u)-\EL(u^*)&=\frac{C_{N,p}}2\int_{\R^N}\frac{\tau_u(z)-\tau_{u^*}(z)}{|z|^N}\dd z\\
&=\frac{C_{N,p}}2\lim_{\varepsilon\downarrow0}D_\varepsilon(u),
\end{aligned}
\end{equation}
where the integral is absolutely convergent. Furthermore, if $0<\varepsilon_1<\varepsilon_2$, then
\begin{equation}\label{4-3}
D_{\varepsilon_1}(u)\ge D_{\varepsilon_2}(u)\ge0.
\end{equation}
Both inequalities in~\eqref{4-3} are strict unless $u$ is a translate of $u^*$.
\end{lemma}
\begin{proof}
By Lemma~\ref{lem3-2}, we have $u^*\in X^p(\R^N)$. Since $\Sp(u^*)=\Sp(u)$, both energies in~\eqref{4-2} are finite. Near the origin,
$$
|\tau_u(z)-\tau_{u^*}(z)|\le\tau_u(z)+\tau_{u^*}(z),
$$
and integrability follows from $u,u^*\in X^p(\R^N)$. For $|z|\ge1$,~\eqref{2-9} and $M_u=M_{u^*}$ yield
$$
|\tau_u(z)-\tau_{u^*}(z)|=|I_{u^*}(z)-I_u(z)|\le I_u(z)+I_{u^*}(z).
$$
The right-hand side is integrable against $|z|^{-N}$ on $|z|\ge1$ by Lemma~\ref{lem2-2}. Subtracting the two identities in~\eqref{2-10} therefore proves the first equality in~\eqref{4-2}.

By~\eqref{3-7}, we have
\begin{equation}\label{4-4}
D_\varepsilon(u)=\int_{\R^N}K_\varepsilon(z)(\tau_u(z)-\tau_{u^*}(z))\dd z.
\end{equation}
Since $K_\varepsilon(z)\le|z|^{-N}$ for $z\ne0$, the dominated convergence theorem proves the second equality in~\eqref{4-2}.

To prove~\eqref{4-3}, let $0<\varepsilon_1<\varepsilon_2$ and set
$$
\ell(r):=(r^2+\varepsilon_1^2)^{-N/2}-(r^2+\varepsilon_2^2)^{-N/2}.
$$
The function $\ell$ is bounded and positive. Moreover, for $r>0$,
$$
\ell'(r)=-Nr\bigl((r^2+\varepsilon_1^2)^{-N/2-1} -(r^2+\varepsilon_2^2)^{-N/2-1}\bigr)<0.
$$
Since $0<\ell\le\kappa_{\varepsilon_1}$, both interactions against $\ell$ are finite. By Proposition~\ref{prop3-1}, we have
$$
D_{\varepsilon_1}(u)-D_{\varepsilon_2}(u)=J_\ell(u^*)-J_\ell(u)\ge0.
$$
Applying the equality assertion of the same proposition to $\ell$ and $\kappa_{\varepsilon_2}$ shows that both inequalities in~\eqref{4-3} are strict unless $u$ is a translate of $u^*$.
\end{proof}

We can now prove Theorem~\ref{thm1-1}. Lemma~\ref{lem4-1} settles the nonnegative case; a comparison with the modulus then handles real-valued functions and determines the sign of the equality cases.
\begin{proof}[Proof of Theorem~\ref{thm1-1}]
The assertion is immediate if $u=0$, so assume that $u\not\equiv0$. Suppose first that $u\ge0$. Lemmas~\ref{lem3-2} and~\ref{lem4-1} yield $u^*\in X^p(\R^N)$ and
$$
\EL(u)-\EL(u^*)=\frac{C_{N,p}}2\lim_{\varepsilon\downarrow0}D_\varepsilon(u)\ge0.
$$
If equality holds, then by~\eqref{4-3}, for every $\varepsilon_0>0$ we have
$$
0\le D_{\varepsilon_0}(u)\le\lim_{\varepsilon\downarrow0}D_\varepsilon(u)=0.
$$
Proposition~\ref{prop3-1} implies that $u=u^*(\cdot-x_0)$ almost everywhere for some $x_0\in\R^N$. The converse follows from translation invariance.

Now let $u$ be real-valued. By the reverse triangle inequality, $|u|\in X^p(\R^N)$, and clearly $\Sp(|u|)=\Sp(u)$. Both energies are finite by Lemma~\ref{lem2-2} and~\eqref{2-3}. Subtracting their expressions in~\eqref{1-4} yields
\begin{equation}\label{4-5}
\begin{aligned}
&\EL(u)-\EL(|u|)\\
&\quad=\frac{C_{N,p}}2\iint_{\R^N\times\R^N} \frac{|u(x)-u(y)|^p-\bigl||u(x)|-|u(y)|\bigr|^p}{|x-y|^N}\dd x\dd y\ge0.
\end{aligned}
\end{equation}
The integrand is nonnegative. For $|x-y|<1$, it is bounded by the short-range integrand of $u$. For $|x-y|\ge1$, integrability follows from~\eqref{2-3} and Lemma~\ref{lem2-2}.

If both $\{u>0\}$ and $\{u<0\}$ have positive measure, the numerator in~\eqref{4-5} is strictly positive on their product. Thus equality in~\eqref{4-5} holds if and only if $u$ does not change sign. Together with the nonnegative case, this implies
$$
\EL(u)\ge\EL(|u|)\ge\EL((|u|)^*)=\EL(u^\#).
$$
If equality holds, $u$ has one sign and $|u|$ is a translate of $u^\#$. Conversely, these properties imply equality by invariance of the energy under translations and multiplication by $-1$.
\end{proof}

The preceding proof also controls the deficit at each fixed regularization scale.
\begin{remark}\label{rem4-2}
Let $u\in X^p(\R^N)$ be nonnegative and satisfy $\Sp(u)<\infty$. By Lemma~\ref{lem4-1}, for every $\varepsilon_0>0$ we have
\begin{equation}\label{4-6}
\EL(u)-\EL(u^*)\ge\frac{C_{N,p}}2 \bigl(J_{\varepsilon_0}(u^*)-J_{\varepsilon_0}(u)\bigr).
\end{equation}
The right-hand side is positive unless $u$ is a translate of $u^*$.
\end{remark}

The assumption $p>1$ is used in the strict supermodularity of $H_p$ in~\eqref{3-1}. At $p=1$, one has $H_1(a,b)=2\min\{a,b\}$, which is not strictly supermodular, so the preceding equality argument does not apply.

\section{Proof of Theorem~\ref{thm1-2}}\label{s5}
By~\cite[Theorem 7.1]{DJS}, the infimum in~\eqref{1-9} is attained by a function $\phi_\Omega\in X_0^p(\Omega)$ satisfying
\begin{equation}\label{5-1}
\norm{\phi_\Omega}_{L^p(\R^N)}=1,\qquad\phi_\Omega>0\quad\text{a.e. in }\Omega.
\end{equation}
These conclusions hold without any regularity assumption on $\partial\Omega$. Since $\phi_\Omega$ vanishes almost everywhere outside the bounded set $\Omega$, we also have $\Sp(\phi_\Omega)<\infty$ by~\eqref{1-6}. Thus Theorem~\ref{thm1-1} applies to $\phi_\Omega$.

\begin{proof}[Proof of Theorem~\ref{thm1-2}]
By translation invariance, we may assume that $B$ is centered at the origin. Let $\phi_\Omega$ be the normalized positive first eigenfunction in~\eqref{5-1}, extended by zero outside $\Omega$. By positivity and equimeasurability, we have
$$
|\{\phi_\Omega>0\}\mathbin\triangle\Omega|=0,\qquad |\{\phi_\Omega^*>0\}\mathbin\triangle B|=0.
$$
Therefore $\phi_\Omega^*\in X_0^p(B)$ and $\norm{\phi_\Omega^*}_{L^p(\R^N)}=1$. By Theorem~\ref{thm1-1},
\begin{equation}\label{5-2}
\lambda_{1,L,p}(B)\le\EL(\phi_\Omega^*)\le\EL(\phi_\Omega)=\lambda_{1,L,p}(\Omega).
\end{equation}
This proves~\eqref{1-10}.

If equality holds, it holds throughout~\eqref{5-2}, and Theorem~\ref{thm1-1} implies that $\phi_\Omega=\phi_\Omega^*(\cdot-x_0)$ almost everywhere for some $x_0\in\R^N$. Taking the positivity sets and using~\eqref{5-1}, we obtain
$$
|\Omega\mathbin\triangle(B+x_0)|=0.
$$
Conversely, if this identity holds, then $X_0^p(\Omega)=X_0^p(B+x_0)$ as spaces of almost-everywhere equivalence classes. Translation invariance yields $\lambda_{1,L,p}(\Omega)=\lambda_{1,L,p}(B)$.

Finally, assume that $\Omega$ has Lipschitz boundary and set $D=B+x_0$. Suppose that $x\in D\setminus\Omega$. Since $D$ is open, a sufficiently small ball centered at $x$ lies in $D$. If $x$ is an interior point of $\Omega^c$, then $|D\setminus\Omega|>0$. If $x\in\partial\Omega$, the same conclusion follows from the local Lipschitz representation of the boundary. Both cases contradict $|\Omega\mathbin\triangle D|=0$, so $D\subset\Omega$. By the same argument, with $\Omega$ and $D$ interchanged and using the smooth boundary of $D$, we obtain $\Omega\subset D$. Hence $\Omega=D$.
\end{proof}

\section{Sharpness of the main theorems}\label{s6}
We conclude the main part of the paper with several observations on sharpness. They concern the logarithmic condition, the assumption $p>1$ in the equality case, and the dependence of the spectral bound on the domain and its volume.

\subsection{The logarithmic finiteness condition}
By Proposition~\ref{prop2-3} and Lemma~\ref{lem3-2}, for every real-valued $u\in X^p(\R^N)$,
\begin{equation}\label{6-1}
\EL(u^\#)\in\R\quad\Longleftrightarrow\quad\Sp(u)<\infty.
\end{equation}
Thus the logarithmic condition in Theorem~\ref{thm1-1} is necessary and sufficient for the rearranged energy to be finite. If it fails, the short-range term of $u^\#$ remains finite, whereas its long-range term equals $-\infty$.

The functions $v_\alpha$ defined in~\eqref{2-15} make the threshold explicit. For every $\alpha>1$, one has $v_\alpha\in W^{1,p}(\R^N)\subset X^p(\R^N)$, while
\begin{equation}\label{6-2}
\Sp(v_\alpha)<\infty\quad\Longleftrightarrow\quad\alpha>2.
\end{equation}
In particular, $W^{1,p}$-regularity alone does not ensure finiteness of the rearranged logarithmic energy.

The power of the logarithm cannot be decreased. Fix $0<\beta<1$ and choose $1+\beta<\alpha\le2$. Since
$$
v_\alpha(x)^p\asymp r^{-N}(\log r)^{-\alpha},\qquad \log^-v_\alpha(x)\asymp\log r\quad\text{as }r=|x|\to\infty,
$$
we have
\begin{equation}\label{6-3}
\int_{\R^N}v_\alpha^p(\log^-v_\alpha)^\beta\dd x<\infty, \qquad\Sp(v_\alpha)=\infty.
\end{equation}
Indeed, in polar coordinates the tail of the first integral is comparable to
$$
\int^\infty\frac{\dd r}{r(\log r)^{\alpha-\beta}},
$$
which is finite because $\alpha-\beta>1$. Hence no smaller positive power of $\log^-|u|$ can replace the logarithm in the finiteness criterion. The same family with $1<\alpha\le2$ shows that $L^p$-integrability alone is insufficient.

For $\alpha>2$, the function $v_\alpha$ satisfies all hypotheses of Theorem~\ref{thm1-1}, although it is positive on $\R^N$ and belongs to no $L^q(\R^N)$ with $0<q<p$. These functions show that neither a restriction on the support nor additional $L^q$-integrability is needed in Theorem~\ref{thm1-1}.

\subsection{The equality case at \texorpdfstring{$p=1$}{p=1}}
The assumption $p>1$ is essential for the equality classification in Theorem~\ref{thm1-1}. To see this, let $\mathcal E_{L,1}$ denote the functional obtained by setting $p=1$ in~\eqref{1-1} and~\eqref{1-4}. Choose $R>r>0$ and $a\ne0$ such that $|a|+r<R$, and set
\begin{equation}\label{6-4}
E:=B_R(0),\qquad F:=B_r(a),\qquad u:=\mathbf1_E+\mathbf1_F.
\end{equation}
Then $F\subset E$ and
$$
u^*=\mathbf1_{B_R(0)}+\mathbf1_{B_r(0)}.
$$
Both functions have finite endpoint energy. Indeed, for a ball $D$,
$$
\int_{\R^N}|\mathbf1_D(x+z)-\mathbf1_D(x)|\dd x =|D\mathbin\triangle(D-z)|\le C_D\min\{|z|,1\}.
$$
Integrating over $|z|<1$ against $|z|^{-N}$ proves short-range integrability. The corresponding term for $u$ is bounded by the sum of those for $\mathbf1_E$ and $\mathbf1_F$. The long-range term is finite because $u$ is bounded and compactly supported, and $H_1(a,b)=2\min\{a,b\}$ for $a,b\ge0$.

For nonnegative $u$, the layer-cake identity reads
$$
|u(x)-u(y)|=\int_0^\infty \bigl|\mathbf1_{\{u>t\}}(x)-\mathbf1_{\{u>t\}}(y)\bigr|\dd t.
$$
In the present example, $\{u>t\}=E$ for $0<t<1$ and $\{u>t\}=F$ for $1<t<2$. The short-range terms therefore add, as do the mass terms and the renormalized long-range integrands. Consequently,
\begin{equation}\label{6-5}
\mathcal E_{L,1}(u) =\mathcal E_{L,1}(\mathbf1_E)+\mathcal E_{L,1}(\mathbf1_F).
\end{equation}
By translation invariance,
\begin{equation}\label{6-6}
\mathcal E_{L,1}(u)=\mathcal E_{L,1}(u^*).
\end{equation}
However, $u$ is not a translate of $u^*$, since its two nontrivial superlevel sets are balls with different centers. Therefore equality does not imply a common center at $p=1$, even for nonnegative compactly supported functions.

\subsection{Equality cases and null modifications of the domain}
Equality in Theorem~\ref{thm1-1} is attained by every admissible nonzero function $\sigma v(\cdot-x_0)$, where $v$ is nonnegative and radially nonincreasing and $\sigma\in\{-1,1\}$. Hence the inequality cannot be improved by a positive additive constant on the whole admissible class. The sign condition is also necessary: by~\eqref{4-5}, $\EL(u)>\EL(|u|)$ whenever both $\{u>0\}$ and $\{u<0\}$ have positive measure.

Balls attain the bound in Theorem~\ref{thm1-2}, so the spectral inequality admits no positive additive improvement at fixed volume. The qualification by null sets in its equality statement is necessary for general open sets. For example, let $B$ be a ball containing the origin and set
\begin{equation}\label{6-7}
\Omega:=B\setminus\{0\}.
\end{equation}
Then $\Omega$ is a bounded open set that is not a ball, but
$$
|\Omega|=|B|,\qquad X_0^p(\Omega)=X_0^p(B),\qquad \lambda_{1,L,p}(\Omega)=\lambda_{1,L,p}(B).
$$
Indeed, the zero-exterior condition in~\eqref{1-8} is imposed almost everywhere. Removing a null set therefore leaves the admissible space unchanged. For domains with Lipschitz boundary, equality modulo null sets implies equality of the domains, as proved in Section~\ref{s5}.

\subsection{Dilation and volume dependence}
To determine the dependence on volume in Theorem~\ref{thm1-2}, we first compute the effect of dilation on the energy. For $t>0$ and $u\in X^p(\R^N)$ with $\Sp(u)<\infty$, set
$$
u_t(x):=t^{-N/p}u(x/t).
$$
With the notation in~\eqref{2-8}, we have $M_{u_t}=M_u$ and $\tau_{u_t}(z)=\tau_u(z/t)$. Dilation preserves $X^p(\R^N)$, since changing the cutoff radius in the defining integral adds or removes only an annulus, on which $\tau_u\le2^pM_u$. It also preserves the logarithmic condition, because
$$
\Sp(u_t)\le\Sp(u)+\frac Np|\log t|M_u.
$$
For $s>0$, define
$$
A_s(u):=\int_{|z|<s}\frac{\tau_u(z)}{|z|^N}\dd z +\int_{|z|\ge s}\frac{\tau_u(z)-2M_u}{|z|^N}\dd z.
$$
Both integrals are finite, and comparison of the cutoff radii yields
$$
A_s(u)=A_1(u)+2M_u|S^{N-1}|\log s.
$$
The change of variables $z=tw$ implies $A_1(u_t)=A_{1/t}(u)$. By~\eqref{2-10} and $C_{N,p}|S^{N-1}|=p$, we obtain
\begin{equation}\label{6-8}
\EL(u_t)=\EL(u)-p\log t\,\norm{u}_{L^p(\R^N)}^p.
\end{equation}
The map $u\mapsto u_t$ is a bijection from $X_0^p(\Omega)$ onto $X_0^p(t\Omega)$ and preserves the $L^p$-norm. Taking infima of the Rayleigh quotients in~\eqref{6-8}, we conclude that
\begin{equation}\label{6-9}
\lambda_{1,L,p}(t\Omega)=\lambda_{1,L,p}(\Omega)-p\log t.
\end{equation}
Writing $\omega_N:=|B_1|$, we may therefore express Theorem~\ref{thm1-2} as
\begin{equation}\label{6-10}
\lambda_{1,L,p}(\Omega)\ge \lambda_{1,L,p}(B_1)-\frac pN\log\!\left(\frac{|\Omega|}{\omega_N}\right).
\end{equation}
Equality holds if and only if $\Omega$ agrees with a ball up to a null set. The coefficient $p/N$ is determined by dilation: for $\Omega=B_R$, formula~\eqref{6-9} yields equality for every $R>0$. Any other coefficient with the same unit-ball normalization would violate the bound for either $R>1$ or $0<R<1$.

The fixed-scale estimate~\eqref{4-6} bounds the energy deficit from below by a regularized Riesz deficit. A bound in terms of the distance to the equality family would require a further quantitative argument.

\begin{remark}\label{rem6-1}
The argument in Section~\ref{s4} and the polarization proof in Appendix~\ref{appA} apply to the same admissible class and establish the same equality statement. In Section~\ref{s4}, monotonicity in Lemma~\ref{lem4-1} reduces equality for the logarithmic energy to equality in a rearrangement inequality at a fixed scale. The strict rearrangement theorem then identifies the center. The same argument proves~\eqref{4-6}.

For the polarization proof, the starting point is an exact deficit for a single reflection. Iterated polarizations lead to the Schwarz comparison, and equality for all affine half-spaces implies radial monotonicity about a common center. To pass to the logarithmic energy, one must also control the long-range interactions. The estimates of Section~\ref{s2} and a cutoff argument ensure this control for functions of arbitrary support.
\end{remark}

\appendix
\section{An alternative proof of Theorem~\ref{thm1-1}}\label{appA}
For completeness, we present a proof of Theorem~\ref{thm1-1} by polarization. We first compute the deficit for a single polarization and characterize its equality cases. We then justify passage to Schwarz rearrangement and show that equality forces the reflection orderings to have a common center.

\subsection{Polarization and the sorting inequality}
Let $H\subset\R^N$ be an open affine half-space, and let $\sigma_H$ denote reflection across $\partial H$. For a real-valued measurable function $u$, define
\begin{equation}\label{A-1}
u^H(x):=\begin{cases}
\max\{u(x),u(\sigma_Hx)\},&x\in H,\\
\min\{u(x),u(\sigma_Hx)\},&x\notin H.
\end{cases}
\end{equation}
The value on $\partial H$ is immaterial. Polarization permutes the two values on each reflection pair and therefore preserves their distribution. In particular,
\begin{equation}\label{A-2}
\norm{u^H}_{L^p(\R^N)}^p=\norm{u}_{L^p(\R^N)}^p, \qquad\Sp(u^H)=\Sp(u).
\end{equation}

For $a,b,c,d\in\R$, set
$$
a_+:=\max\{a,b\},\quad a_-:=\min\{a,b\},\quad c_+:=\max\{c,d\},\quad c_-:=\min\{c,d\},
$$
and define
\begin{equation}\label{A-3}
\begin{aligned}
Q_p(a,b,c,d):={}&|a-c|^p+|b-d|^p\\
&-|a_+-c_+|^p-|a_--c_-|^p.
\end{aligned}
\end{equation}
For $p>1$, we have
\begin{equation}\label{A-4}
Q_p(a,b,c,d)\ge0,\qquad Q_p(a,b,c,d)>0\quad\Longleftrightarrow\quad(a-b)(c-d)<0.
\end{equation}
Indeed, if $(a-b)(c-d)\ge0$, the two sums in~\eqref{A-3} coincide. If $a>b$ and $d>c$, then
\begin{align*}
Q_p(a,b,c,d) &=|a-c|^p+|b-d|^p-|a-d|^p-|b-c|^p\\
&=p(p-1)\int_b^a\int_c^d|s-t|^{p-2}\dd t\dd s>0.
\end{align*}
The other case follows by interchanging the two entries in each pair. When $1<p<2$, the possible singularity at $s=t$ is locally integrable, as in~\eqref{3-1}.

\subsection{The inequality and equality for one polarization}
We next apply the sorting inequality to the short-range form and to the full logarithmic energy. The resulting deficit identity also determines the equality cases for a single polarization.
\begin{proposition}\label{propA-1}
Let $1<p<\infty$, let $u\in X^p(\R^N)$ be real-valued, and let $H$ be an open affine half-space. Then $u^H\in X^p(\R^N)$ and
\begin{equation}\label{A-5}
[u^H]_{X^p}^p\le[u]_{X^p}^p.
\end{equation}
If $\Sp(u)<\infty$, both logarithmic energies are finite and
\begin{equation}\label{A-6}
\EL(u^H)\le\EL(u).
\end{equation}
Equality holds in~\eqref{A-6} if and only if
\begin{equation}\label{A-7}
u(x)\ge u(\sigma_Hx)\quad\text{for a.e. }x\in H,
\end{equation}
or the reverse inequality holds for almost every $x\in H$. Equivalently, equality holds if and only if $u^H=u$ almost everywhere or $u^H=u\circ\sigma_H$ almost everywhere.
\end{proposition}
\begin{proof}
For $x,y\in H$, write
$$
a=u(x),\quad b=u(\sigma_Hx),\qquad c=u(y),\quad d=u(\sigma_Hy).
$$
Reflection geometry implies $|x-y|<|x-\sigma_Hy|$. Let $k$ be a bounded, integrable, radial nonincreasing kernel. Grouping the interaction into the four pairs associated with $x,y\in H$, we obtain
\begin{equation}\label{A-8}
\begin{aligned}
&\iint_{\R^N\times\R^N}k(x-y) \bigl(|u(x)-u(y)|^p-|u^H(x)-u^H(y)|^p\bigr)\dd x\dd y\\
&\quad=\iint_{H\times H} \bigl(k(x-y)-k(x-\sigma_Hy)\bigr)Q_p(a,b,c,d)\dd x\dd y\ge0.
\end{aligned}
\end{equation}
This identity follows because sorting preserves the sum of the four differences. Apply~\eqref{A-8} to
$$
k_\delta(z):=\mathbf1_{\{|z|<1\}}\min\{\delta^{-N},|z|^{-N}\},\qquad0<\delta<1.
$$
Letting $\delta\downarrow0$, the monotone convergence theorem proves~\eqref{A-5}. By~\eqref{A-2}, Lemma~\ref{lem2-2} and~\eqref{2-3}, both energies are finite when $\Sp(u)<\infty$.

We may now group the logarithmic energies into the same reflection pairs. The mass terms cancel. The sum of the long-range subtraction terms on each group is
$$
-(f_1+f_2)(|a|^p+|b|^p+|c|^p+|d|^p),
$$
where
$$
f_1:=\frac{\mathbf1_{\{|x-y|\ge1\}}}{|x-y|^N},\qquad f_2:=\frac{\mathbf1_{\{|x-\sigma_Hy|\ge1\}}}{|x-\sigma_Hy|^N}.
$$
This expression is unchanged by sorting. Consequently,
\begin{equation}\label{A-9}
\begin{aligned}
\EL(u)-\EL(u^H) =\frac{C_{N,p}}2\iint_{H\times H} &\left(\frac1{|x-y|^N}-\frac1{|x-\sigma_Hy|^N}\right)\\
&\times Q_p(a,b,c,d)\dd x\dd y\ge0.
\end{aligned}
\end{equation}
These regroupings are justified by finiteness of the short-range terms and absolute convergence of the long-range terms. The kernel difference in~\eqref{A-9} is strictly positive almost everywhere on $H\times H$.

Set
$$
A:=\{x\in H:u(x)>u(\sigma_Hx)\},\qquad B:=\{x\in H:u(x)<u(\sigma_Hx)\}.
$$
If both sets have positive measure,~\eqref{A-4} implies that the integrand in~\eqref{A-9} is strictly positive on $A\times B$. Equality therefore forces $|A|=0$ or $|B|=0$, which is precisely the asserted ordering. Conversely, either ordering makes $u^H$ equal to $u$ or its reflection. Reflection invariance of the energy then proves equality.
\end{proof}

\begin{remark}\label{remA-2}
For real-valued $u$ under the hypotheses of Theorem~\ref{thm1-1}, Proposition~\ref{propA-1} and~\eqref{4-5} yield
\begin{equation}\label{A-10}
\EL((|u|)^H)\le\EL(|u|)\le\EL(u).
\end{equation}
If $u\not\equiv0$, equality holds throughout if and only if $u$ has one sign and $|u|$ satisfies one of the two orderings across $\partial H$ in Proposition~\ref{propA-1}.

Equality for signed polarization~\eqref{A-1} does not require one sign, since every reflection-invariant function is fixed by that operation. Nor does equality for a single polarization imply radial symmetry. For example, every nonnegative smooth compactly supported function whose support is contained in $H$ satisfies $u^H=u$, regardless of the shape of its level sets.
\end{remark}

\subsection{Passage to Schwarz rearrangement}
For nonnegative $u\in L^p(\R^N)$, the approximation theorem of Brock and Solynin~\cite{BS} and the constructive version of Hajaiej~\cite{HajPol} provide finite iterated polarizations $u_n$ toward half-spaces containing the origin such that
\begin{equation}\label{A-11}
u_n\longrightarrow u^*\quad\text{in }L^p(\R^N).
\end{equation}
More precisely, the constructive scheme uses a dense sequence $(H_j)$ of half-spaces containing the origin in their interior, with $u_0=u$ and $u_{n+1}=u_n^{H_1\cdots H_{n+1}}$. The repeated finite blocks are part of the convergence statement; an arbitrary sequence of half-spaces need not converge to $u^*$.

To use~\eqref{A-11}, we first consider functions with bounded support, for which the long-range term is continuous under $L^p$-convergence. We then use the estimates of Section~\ref{s2} to remove the support restriction.
\begin{proposition}\label{propA-3}
Let $u\in X^p(\R^N)$ be real-valued and satisfy $\Sp(u)<\infty$. Then $u^\#\in X^p(\R^N)$, $\Sp(u^\#)=\Sp(u)$, and $\EL(u^\#)\le\EL(u)$.
\end{proposition}
\begin{proof}
Suppose first that $u\ge0$ and $u=0$ outside $B_R$. All iterates in~\eqref{A-11} also vanish outside $B_R$, since polarization toward a half-space containing the origin places the larger value at the point closer to the origin. By Proposition~\ref{propA-1},
\begin{equation}\label{A-12}
\EL(u_n)\le\EL(u)\quad\text{for every }n.
\end{equation}
On functions supported in a common bounded set, the long-range term is continuous under $L^p$-convergence. Indeed, its numerator vanishes when either argument is zero, so the interaction is confined to a bounded product set. The kernel is bounded on $|x-y|\ge1$, and the standard $L^p$ estimates for differences of $p$th powers imply continuity of the integral. The mass terms also converge. After passing to an almost everywhere convergent subsequence, Fatou's lemma for the short-range term yields $u^*\in X^p(\R^N)$ and
\begin{equation}\label{A-13}
\EL(u^*)\le\liminf_{n\to\infty}\EL(u_n)\le\EL(u).
\end{equation}

Now let $u\ge0$ have arbitrary support and satisfy $\Sp(u)<\infty$. Choose a smooth radial nonincreasing function $\chi$, equal to $1$ on $B_1$ and to $0$ outside $B_2$, with $0\le\chi\le1$. Set $\chi_R(x):=\chi(x/R)$ and $u_R:=\chi_Ru$ for $R\ge1$. Then $u_R\uparrow u$ almost everywhere and $u_R\to u$ in $L^p(\R^N)$. To verify convergence in $X^p(\R^N)$, put $w_R:=(1-\chi_R)u$. By the Lipschitz bound for $\chi_R$,
\begin{align*}
[w_R]_{X^p}^p &\le2^{p-1}\iint_{|x-y|<1} \frac{(1-\chi_R(x))^p|u(x)-u(y)|^p}{|x-y|^N}\dd x\dd y\\
&\quad+CR^{-p}\norm{u}_{L^p(\R^N)}^p,
\end{align*}
where $C$ is independent of $R$. The first term tends to zero by dominated convergence, so $u_R\to u$ in $X^p(\R^N)$.

The function $H_p$ is increasing in each nonnegative argument. Since $u_R\uparrow u$ and $u_R^*\uparrow u^*$ almost everywhere, monotone convergence applies to both long-range interactions. Their limits are finite by Lemma~\ref{lem2-2}, since $\Sp(u^*)=\Sp(u)$. It follows that $\EL(u_R)\to\EL(u)$, while Fatou's lemma for the rearranged short-range term yields $u^*\in X^p(\R^N)$ and
\begin{equation}\label{A-14}
\EL(u^*)\le\liminf_{R\to\infty}\EL(u_R^*) \le\lim_{R\to\infty}\EL(u_R)=\EL(u).
\end{equation}
Here the second inequality follows from the bounded-support case applied to $u_R$.

For real-valued $u$, apply the nonnegative case to $|u|$ and use~\eqref{4-5}. Equimeasurability implies $\Sp(u^\#)=\Sp(u)$.
\end{proof}

\subsection{Equality and the common center}
It remains to identify the equality cases. By Proposition~\ref{propA-1}, equality for every affine half-space imposes a reflection ordering across each hyperplane. The next lemma shows that these orderings determine a common center.
\begin{lemma}\label{lemA-4}
Let $0\ne v\in L^p(\R^N)$ be nonnegative. Suppose that, for every open affine half-space $H$, either $v(x)\ge v(\sigma_Hx)$ almost everywhere on $H$, or the reverse inequality holds almost everywhere on $H$. Then $v=v^*(\cdot-x_0)$ almost everywhere for some $x_0\in\R^N$.
\end{lemma}
\begin{proof}
Set $w:=v^p\in L^1(\R^N)$, $\Phi(x):=e^{-|x|^2}$ and $F:=\Phi*w$. The function $F$ is positive and continuous and tends to zero at infinity. Hence it attains its maximum at some point $x_0\in\R^N$. The reflection orderings of $v$ also hold for $w$.

Let $H$ contain $x_0$ in its interior. If $w(y)\le w(\sigma_Hy)$ almost everywhere on $H$, with strict inequality on a set of positive measure, then
\begin{align*}
F(x_0)-F(\sigma_Hx_0) =\int_H&\bigl(\Phi(x_0-y)-\Phi(x_0-\sigma_Hy)\bigr)\\
&\times\bigl(w(y)-w(\sigma_Hy)\bigr)\dd y<0.
\end{align*}
Indeed, the first factor is strictly positive for $y\in H$. This contradicts maximality at $x_0$. Therefore $v(x)\ge v(\sigma_Hx)$ almost everywhere on $H$ for every half-space containing $x_0$ in its interior. Thus $v^H=v$ for all such half-spaces. By the characterization of Schwarz rearrangement through inward polarizations~\cite{HajPol}, translated from the origin to $x_0$, $v$ is radially nonincreasing about $x_0$. Equimeasurability identifies its profile with that of $v^*$.
\end{proof}

We now combine this geometric fact with the polarization inequality to recover the equality statement of Theorem~\ref{thm1-1}.
\begin{proposition}\label{propA-5}
Under the hypotheses of Theorem~\ref{thm1-1}, equality holds for $u\not\equiv0$ if and only if there exist $x_0\in\R^N$ and $\sigma\in\{-1,1\}$ such that $u=\sigma u^\#(\cdot-x_0)$ almost everywhere.
\end{proposition}
\begin{proof}
Suppose first that $u\ge0$ and $\EL(u)=\EL(u^*)$. For every open affine half-space $H$, Propositions~\ref{propA-1} and~\ref{propA-3} yield
$$
\EL(u^*)\le\EL(u^H)\le\EL(u),
$$
since $(u^H)^*=u^*$. Equality therefore holds for every polarization. By Proposition~\ref{propA-1}, $u$ satisfies one of the two reflection orderings for each $H$. Lemma~\ref{lemA-4} implies that $u=u^*(\cdot-x_0)$ almost everywhere for some $x_0\in\R^N$. The converse follows from translation invariance.

For real-valued $u$, equality in
$$
\EL(u)\ge\EL(|u|)\ge\EL(u^\#)
$$
forces $u$ to have one sign by~\eqref{4-5}. The nonnegative case then shows that $|u|$ is a translate of $u^\#$. Conversely, these properties imply equality by invariance of the energy under translations and multiplication by $-1$.
\end{proof}

\subsection{Comparison of the two arguments}
Equality for one polarization requires only an ordering across the corresponding hyperplane. It allows nonradial functions and, for signed polarization, sign-changing functions. In the proof above, radial monotonicity follows from equality for every affine half-space and Lemma~\ref{lemA-4}. The sign condition follows separately from~\eqref{4-5}.

The passage from iterated polarizations to the logarithmic energy also requires more than $L^p$-convergence. For nonnegative functions, the long-range contribution is negative, so Fatou's lemma for the short-range term alone does not imply lower semicontinuity of the full energy. Bounded support ensures continuity of the long-range term in~\eqref{A-13}; the estimates of Section~\ref{s2} and the cutoff argument in~\eqref{A-14} remove this restriction. Thus the polarization proof applies under the same logarithmic condition as Theorem~\ref{thm1-1}, without additional fractional regularity or $L^q$-integrability for $q<p$. In particular, it includes the functions $v_\alpha$ with $\alpha>2$ from Remark~\ref{rem2-4}. In Section~\ref{s4}, the same equality conclusion follows from monotonicity of the regularized interaction deficits, and the proof also establishes~\eqref{4-6}.

Finally, $p>1$ is essential for strictness in~\eqref{A-4}. At $p=1$, opposite orderings need not produce a positive deficit. The nested-ball example in~\eqref{6-4} and~\eqref{6-6} shows that the common-center equality classification fails at that endpoint.

\subsection*{Data availability}
No data were used in the research described in this article.
\subsection*{Competing interests}
The authors declare that they have no competing interests.

\end{document}